\documentclass[12pt]{amsart}
\usepackage{amsmath,amssymb}
\usepackage{amsfonts}
\usepackage{amsthm}
\usepackage{latexsym}
\usepackage{graphicx}

\def\C{\mathbb{C}}

\def\vv<#1>{\langle#1\rangle}

\def\XXint#1#2{\setbox0=\hbox{$#1{#2}{\int}$}{#2}\kern-.5\wd0 }

\def\XXint#1#2#3{{\setbox0=\hbox{$#1{#2#3}{\int}$}
     \vcenter{\hbox{$#2#3$}}\kern-.5\wd0}}

\def\vv<#1>{\langle#1\rangle}
\newtheorem{thm}{Theorem}[section]

\newtheorem{lem}{Lemma}[section]

\theoremstyle{definition}

\theoremstyle{remark}

\newtheorem{rem}{Remark}[section]

\numberwithin{equation}{section}

\begin{document}

\title{A Liouville property of holomorphic maps}


\begin{abstract} In this article, we prove a Liouville property of
holomorphic maps from a complete K\"ahler manifold with nonnegative
holomorphic  bisectional curvature to a complete simply connected
K\"ahler manifold with a certain assumption on the sectional
curvature.
\end{abstract}

\renewcommand{\subjclassname}{\textup{2000} Mathematics Subject Classification}
\author{ Chengjie Yu}
\address{Department of Mathematics, Shantou University, Shantou, Guangdong, China } \email{cjyu@stu.edu.cn}

\date{Jan 2010}

\maketitle \markboth{ Chengjie Yu}
 {Liouville property}

 Liouville property is an interesting topic in analysis since Liouville found that a bounded holomorphic function on the
 complex plane must be a constant function. It has been studied by many geometers. For example, Yau \cite{Yau} studied Liouville
 property of harmonic functions on complete Riemannian manifold with nonnegative Ricci curvature and Tam \cite{Tam} studied
 Liouville property of harmonic maps. By Yau's Schwartz lemma(\cite{Ya2}), every holomorphic map from a complete K\"ahler manifold
 with nonnegative Ricci curvature to a complete K\"ahler manifold with holomorphic bisectional curvature not greater than a negative
 constant is a constant map. This is a Liouville property of holomorphic maps.

 In this article,  we prove the following Liouville
 property of holomorphic maps.
 \begin{thm}
 Let $M$ be a complete K\"ahler manifold with nonnegative
 bisectional curvature and $N$ be a complete simply connected K\"ahler
 manifold with sectional curvature $\leq -\frac{c}{(1+r)^2}$  where
 $c$ is some positive constant and $r$ is  the distance function of $N$ with respect to a fixed point. Then,
 any holomorphic map form $M$ to $N$ must be constant.
 \end{thm}

 In order to prove the theorem, we need the following lemma on the existence of plenty of bounded strictly plurisubharmonic
 functions so that they separate points on the manifold.
 \begin{lem}
 Let $M$ be a complete simply connected K\"ahler manifold with sectional curvature $\leq -\frac{c}{(1+r)^2}$  where
 $c$ is  some positive constant and $r$ is  the distance function of $M$ with respect to a fixed point. Then, for any point $p\in M$, there is
 a bounded continuous strictly plurisubharmonic function $\phi$ such that $\phi>0$ on $M\setminus\{p\}$ and $\phi(p)=0$.
 \end{lem}
 \begin{proof} We construct the strictly plurisubharmonic function with the form $\phi=f(r)$ where $r$ is the distance function
 to $p$ and $f$ is a function to be determined.

 By Hessian comparison, we have
 \begin{equation}\label{eqn-hc}
 D^2r\geq\frac{\sqrt c}{1+r}\coth\bigg(\frac{\sqrt cr}{1+r}\bigg)(g-dr\otimes dr).
 \end{equation}
 Let $e_1=\frac{1}{\sqrt{2}}(\nabla r-J\nabla r )$ and $e_1,e_2,\cdots,e_n$
 be a parallel unitary frame along geodesic rays emanating from $p$. Then
 \begin{equation}
 r_{1}=\frac{1}{\sqrt 2}
 \end{equation}
 and $r_{\alpha}=0$ for any $\alpha>1$.

 By \eqref{eqn-hc}, when assuming $f'>0$, we have
 \begin{equation}
 \begin{split}
 \phi_{\alpha\bar\beta}=&f'(r)r_{\alpha\bar\beta}+f''(r)r_\alpha r_{\bar \beta}\\
 \geq&f'(r)\frac{\sqrt c}{1+r}\coth\bigg(\frac{\sqrt cr}{1+r}\bigg)(g_{\alpha\bar\beta}-r_\alpha r_{\bar\beta})\\
 &+f''(r)r_\alpha r_{\bar \beta}\\
 :=&\psi_{\alpha\bar\beta}
 \end{split}
 \end{equation}
 Note that when $\alpha\neq \beta$, we have $\psi_{\alpha\bar \beta}=0$, when $\alpha>1$,
 \begin{equation}
 \psi_{\alpha\bar\alpha}=\frac{\sqrt c}{1+r}\coth\bigg(\frac{\sqrt cr}{1+r}\bigg)f'(r),
 \end{equation}
 and
 \begin{equation}
 \psi_{1\bar 1}=\frac{1}{2}\bigg[\frac{\sqrt c}{1+r}\coth\bigg(\frac{\sqrt cr}{1+r}\bigg)f'(r)+f''(r)\bigg]
 \end{equation}
 So, if $f'(r)>0$ and
 \begin{equation}\label{eqn-f-2}
 [\log
 f']'(r)>-\frac{\sqrt c}{1+r}\coth\bigg(\frac{\sqrt cr}{1+r}\bigg),
 \end{equation}
 then $\phi$ is strictly plurisubharmonic. Note that
 \begin{equation}
\sqrt c\coth\bigg(\frac{\sqrt cr}{1+r}\bigg)\to \sqrt c \coth \sqrt
c
 \end{equation}
 as $r\to \infty$ and $\sqrt c\coth \sqrt c>1$. So, there are two positive numbers $R$ and $\delta$, such that
\begin{equation}
\frac{\sqrt c}{1+r}\coth\bigg(\frac{\sqrt cr}{1+r}\bigg)>\frac{1+\delta}{1+r}
\end{equation}
for any $r\geq R$.

Let $h(r)$ be a continuous function on $[0,\infty)$ such that
\begin{equation}
\left\{\begin{array}{ll}h(r)=\frac{1+\delta}{1+r}<\frac{\sqrt c}{1+r}\coth\bigg(\frac{\sqrt cr}{1+r}\bigg)& r\geq R\\ h(r)<\frac{\sqrt c}{1+r}\coth\bigg(\frac{\sqrt cr}{1+r}\bigg)& \mbox{for}\ r\leq R  \end{array}\right.
\end{equation}
Then
\begin{equation}
f(r)=\int_0^{r}e^{-\int_0^{s}h(t)dt}ds
\end{equation}
satisfies \eqref{eqn-f-2} and $f'>0$. So $\phi=f(r)$ is a strictly plurisubharmonic function. Moreover,
\begin{equation}
\begin{split}
f(r)=&\int_0^{r}e^{-\int_0^{R}h(t)dt-\int_R^{s}h(t)dt}ds\\
=&\int_0^{r}e^{-\int_0^{R}h(t)dt-\int_R^{s}\frac{1+\delta}{1+t}dt}ds\\
=&(1+R)^{1+\delta}e^{-\int_0^{R}h(t)dt}\int_0^{s}\frac{1}{(1+s)^{1+\delta}}ds\\
\leq&(1+R)^{1+\delta}e^{-\int_0^{R}h(t)dt}\int_0^{\infty}\frac{1}{(1+s)^{1+\delta}}ds.
\end{split}
\end{equation}
So, $f$ is bounded. It is clear that $f(0)=0$ and $f(r)>0$ for any $r>0$. Therefore, $\phi=f(r)$ is
a bounded continuous strictly plurisubharmonic function on $M$ with $\phi>0$ on $M\setminus \{p\}$ and $\phi(p)=0$.
 \end{proof}
\begin{rem}
Similar results can be found in \cite{GW}. Our treatment here is different with that in \cite{GW}.
\end{rem}

 \begin{proof}[Proof of the theorem]
 We proceed by contradiction. Let $f:M\to N$ be a nonconstant
 holomorphic map. Then, there are two points $p,q$ in $M$, such that
 $f(p)\neq f(q)$. By the lemma above, there is a
 continuous plurisubharmonic function $\phi$ on $N$ such that $\phi(f(p))=0$
 , $\phi(f(q))>0$. Then,
 $\phi(f)$ is a continuous bounded plurisubharmonic function on $M$ that is
 nonconstant. This contradicts Theorem 3.2 in \cite{NT1}.
 \end{proof}
 \begin{rem}
The unitary invariant metric on $\C^n$ constructed by
 Seshadri \cite{Seshadri2006} has negative sectional curvature and sectional curvature $\leq -\frac{1}{r^2\log
 r}$ outside a compact subset. This means that the curvature decay rate can not be raised to be
 greater than $2$.
 \end{rem}


\begin{thebibliography}{99}
\bibitem{GW}Greene, R. E.; Wu, H. {\sl Function theory on manifolds which possess a pole.} Lecture Notes in Mathematics, 699. Springer, Berlin, 1979. ii+215 pp. ISBN: 3-540-09108-4.
\bibitem{NT1}Ni, Lei; Tam, Luen-Fai. {\sl Plurisubharmonic functions and the structure of complete K\"ahler manifolds with nonnegative curvature.} J. Differential Geom. 64 (2003), no. 3, 457--524.
\bibitem{Seshadri2006}   Seshadri, H., {\sl Negative sectional curvature and the product complex structure}, Math. Res. Lett.  \textbf{13} (2006),  495-500.
\bibitem{Tam}Tam, Luen-Fai.{\sl Liouville properties of harmonic maps.} Math. Res. Lett. 2 (1995), no. 6, 719--735.
\bibitem{Yau}Yau, Shing Tung.{\sl Harmonic functions on complete Riemannian manifolds.} Comm. Pure Appl. Math. 28 (1975), 201--228.
\bibitem{Ya2}Yau, Shing Tung {\sl A general Schwarz lemma for K\"ahler manifolds.} Amer. J. Math. 100(1978), no. 1, 197--203.
\end{thebibliography}
\end{document}